\documentclass[12pt] {article}
\usepackage{amsmath,amsthm}
\usepackage{amssymb,latexsym}

\title{On extendability by continuity of valuations on convex polytopes.}
\date{}
\author{Semyon Alesker \footnote{Partially supported by ISF grant 701/08 and 1447/12.}
\\  { \normalsize Department of Mathematics, Tel Aviv University, Ramat Aviv}
\\  { \normalsize 69978 Tel Aviv, Israel }
\\ {\normalsize e-mail: semyon@post.tau.ac.il}}

\def\RR{\mathbb{R}}
\def\CC{\mathbb{C}}

\def\ZZ{\mathbb{Z}}

\def\PP{\mathbb{P}}

\def\eps{\varepsilon}

\def\ome{\omega}
\def\Ome{\Omega}
\def\lam{\lambda}

\def\to{\longrightarrow}

\def\qed { Q.E.D. }

\swapnumbers
\newtheorem{theorem}{Theorem}[section]

\newtheorem{lemma}[theorem]{Lemma}
\newtheorem{proposition}[theorem]{Proposition}
\newtheorem{claim}[theorem]{Claim}
\theoremstyle{definition}
\newtheorem{example}[theorem]{Example}
\newtheorem{definition}[theorem]{Definition}
\newtheorem{remark}[theorem]{Remark}
 \def\cb{{\cal B}} 
\def\cd{{\cal D}}  \def\cf{{\cal F}}
\def\cg{{\cal G}} \def\ch{{\cal H}} \def\ci{{\cal I}}
 \def\ck{{\cal K}} \def\cl{{\cal L}}
\def\cm{{\cal M}}  \def\co{{\cal O}}
\def\cp{{\cal P}} 

 \def\ct{{\cal T}} \def\cu{{\cal U}}

\def\dl{{\cal D}_{\lambda}}
\def\cm{{\cal M}}
\def\pt{\partial}

\numberwithin{equation}{section}
\begin{document}
\maketitle

\begin{abstract}
There is a well known construction of weakly continuous valuations
on convex compact polytopes in $\RR^n$. In this paper we investigate
when a special case of this construction gives a valuation which
extends by continuity in the Hausdorff metric to all convex compact
subsets of $\RR^n$. It is shown that there is a necessary condition
on the initial data for such an extension. In the case of $\RR^3$
more explicit results are obtained.
\end{abstract}

\tableofcontents

\def\gk{G_k(V)}
\setcounter{section}{-1}
\section{Introduction.} The aim of this paper is to
discuss a well known construction of valuations on convex compact
polytopes in $\RR^n$ and study the question when the constructed
valuations can be extended by continuity in the Hausdorff metric
from polytopes to the class of all convex compact sets. We show that
in general there is a non-trivial obstruction to such an extension.
The case of 3-dimensional space is studied in greater detail, and we
obtain necessary and sufficient conditions on the initial data of the
construction under which the valuations do extend by continuity to
all convex compact sets. The main results of the paper are Theorems
\ref{T:main-res1}, \ref{T:main-res2}, \ref{T:dim-3-main} below.

Let $V$ be a finite dimensional real vector space, $n=\dim V$. Let
$\cp(V)$ denote the family of all convex compact polytopes in $V$,
and let $\ck(V)$ denote the family of all non-empty convex compact subsets in
$V$.
\begin{definition}
A valuation on $\cp(V)$ (resp. $\ck(V)$) is a functional
$$\phi:\cp(V)\to \CC$$
$$(\mbox{ resp. } \phi:\ck(V)\to \CC)$$
which is additive in the following sense: for any $A, B\in \cp(V)$
(resp. $\ck(V)$) such that the union $A\cup B$ is also convex, one
has
$$\phi(A\cup B)=\phi(A)+\phi(B)-\phi(A\cap B).$$
\end{definition}

As a generalization of the notion of measure, valuations on convex
sets have long played an important role in geometry. Since the
middle 1990s' we have witnessed a breakthrough in the structure
theory of valuations \cite{klain1},
\cite{schneider-96},\cite{alesker-gafa}-\cite{alesker-survey},
\cite{bernig-broecker}-\cite{bernig-fu-convolution},
\cite{ludwig}-\cite{ludwig-reitzner-2},\cite{schuster-duke}. This
progress in turn has led to immense advances in the integral
geometry of isotropic (in particular, Hermitian) spaces (see
\cite{alesker-jdg-03},\cite{alesker-IG},\cite{bernig-fu-ann},\cite{bernig-fu-solanes})
and an understanding why certain classical notions from geometric
analysis are indeed fundamental (see
\cite{haberl-parapatits},\cite{ludwig-reitzner},\cite{ludwig},\cite{ludwig-fischer},\cite{schuster-adv-08}).

\hfill

\begin{definition}
A valuation $\phi$ on $\ck(V)$ is called {\itshape continuous} if
$\phi$ is continuous in the Hausdorff metric on $\ck(V)$.
\end{definition}
The notion of {\itshape weak continuity} of a valuation is often
more appropriate for polytopes; it is due to Hadwiger
\cite{hadwiger-52}. Let us recall it. Fix $\xi:=\{\xi_1,\dots,
\xi_s\}\subset V^*$ an $s$-tuple of linear functionals on $V$. Let
$\cp_\xi$ denote the family of polytopes in $\cp(V)$ of the form
$$P_\xi(y)=\{x\in V|\, \xi_i(x)\leq y_i \forall i=1,\dots, s\}$$
where $y=(y_1,\dots,y_s)\in \RR^s$ are such that $P_\xi(y)$ are 
non-empty and compact.
\begin{definition}
A valuation $\phi$ on $\cp(V)$ is called {\itshape weakly
continuous} if for any $\xi$ the function
$$y\mapsto \phi(P_\xi(y))$$
is continuous with respect to $y$ in the set where $P_\xi(y)$ is
non-empty and compact.
\end{definition}

Clearly the restriction to $\cp(V)$ of a continuous valuation on
$\ck(V)$ is weakly continuous. The first question studied in this
article is the converse one: when does a weakly continuous valuation
on $\cp(V)$ extend to a continuous valuation on $\ck(V)$? Clearly
this extension is unique if it exists.

Weakly continuous translation invariant valuations on $\cp(V)$
were described completely by P. McMullen \cite{mcmullen-83}. All
such valuations are obtained by a general relatively explicit
construction. Let us describe a particular case of it which will
be studied in this article in greater detail.

\hfill


Let us fix a Euclidean metric on $V$. Let us denote by $G_k(V)$ the
manifold of pairs $(E,l)$ where $E\subset V$ is a $k$-dimensional
linear subspace of $V$, and $l$ is a line orthogonal to $E$.


Let $f:\gk\to \CC$ be a continuous function.
For a polytope $P\in \cp(V)$ let us denote by $\cf_k(P)$ the set of
$k$-faces of $P$. For any $k$-face $F\in \cf_k(P)$ let us denote by
$\bar F$ the (only) $k$-dimensional linear subspace of $V$
containing a translate of $F$. Let us denote by $\gamma_F$
the exterior angle of $P$ at $F$, that is $\gamma_F$ is the subset
of the unit sphere $S(\bar F^\perp)$ of $\bar F^\perp$ consisting of
all $n\in S(\bar F^\perp)$ such that for some (equivalently, any)
$x$ belonging to the relative interior of $F$ the scalar product
$(n,p-x)$ is non-positive for any $p\in P$.

Finally for $P\in \cp(V)$ define
\begin{eqnarray}\label{E:expression}
\phi_f(P)=\sum_{F\in \cf_k(P)}vol (F) \int_{\gamma_F}f(\bar F,l)dl,
\end{eqnarray}
where $dl$ denotes the Lebesgue measure on the sphere $S(\bar
F^\perp)$ such that the total measure of $S(\bar F^\perp)$ is equal
to 1 (here we consider $f(\bar F,\cdot)$ as an even function on the
unit sphere $S(\bar F^\perp)$).

By the (easy part of) P. McMullen's theorem \cite{mcmullen-83},
$\phi_f$ is a translation invariant even weakly continuous valuation
on $\cp(V)$.

\begin{example}\label{example1}
Take $f\equiv 1$. Then $\phi_f$ is proportional to the $k$-th
intrinsic volume $V_k$ (see the book \cite{schneider-book} for this
notion).
\end{example}

Let us mention that expressions similar to (\ref{E:expression}) have
been studied in recent preprints \cite{hinderer-hug-weil},
\cite{schneider-preprint}.

\hfill

The goal of this paper is to study the following question: when does
the valuation $\phi_f$ extend by continuity in the Hausdorff metric
to $\ck(V)$? It turns out that such an extension does not always
exist. Notice that by a result of Groemer \cite{groemer} if $\phi_f$
extends by continuity to $\ck(V)$ then this extension will be
automatically a valuation on $\ck(V)$. In this case we will denote
this extension by the same symbol $\phi_f$.

Our first main result is as follows.
\begin{theorem}\label{T:main-res1}
Let $n:=\dim V\geq 3$ and $1\leq k\leq n-2$. There exists a real
analytic function $f\colon G_k(V)\to \CC$ such that $\phi_f$ does
\underline{not} extend to $\ck(V)$ by continuity.
\end{theorem}
\begin{remark}
(1) It is well known that for $k=n-1$ the valuation $\phi_f$ does
extend by continuity to $\ck(V)$ for any continuous function $f$.

(2) Obviously the simplest case when the assumptions of Theorem
\ref{T:main-res1} are satisfied is $n=3,k=1$. It turns out that in
this case one has a much more precise statement, see Theorem
\ref{T:dim-3-main} below.

(3) We would like to emphasize the following phenomenon: the class
of smoothness of the function $f$ is not really important, it can be
real analytic, and in fact even $O(n)$-finite.\footnote{Recall the
definition of $O(n)$-finite vectors in a representation of $O(n)$ in
a vector space $\cu$: a vector $\xi\in\cu$ is called $O(n)$-finite
if it is contained in a finite dimensional $O(n)$-invariant linear
subspace. It is well known that if $\cu$ is the space of (smooth,
continuous, or $L^2$) functions on a homogeneous space of $O(n)$,
then all $O(n)$-finite functions are real analytic.}

(4) By Proposition \ref{P:closedness} below, the subset of the space
$C^{\infty}(G_k(V))$ of infinitely smooth functions $f$ such that
$\phi_f$ admits an extension to $\ck(V)$ by continuity is a
{\itshape closed} linear subspace.

(5) In a very recent preprint \cite{hinderer-hug-weil} Theorem
\ref{T:main-res1} was proved by a very different method.
\end{remark}

Now we are going to explain some more precise forms of Theorem
\ref{T:main-res1}. In the rest of the paper we will always assume
that $f$ is infinitely smooth.

Let us denote by $C^{\infty}(X)$ the space of infinitely smooth
functions on a manifold $X$. Let us denote by $Gr_k(V)$ the
Grassmannian of $k$-dimensional linear subspaces of $V$. Let us
define the linear map
$$S\colon C^{\infty}(G_k(V))\to C^\infty(Gr_k(V))$$
by $$S(g)(E)=\int_{l\in \PP(E^\perp)}g(E,l)dl,$$ where
$\PP(E^\perp)$ denotes the set of lines in $E^\perp$ (i.e. the
projectivization of $E^\perp$), $dl$ is the normalized Lebesgue
measure on $\PP(E^\perp)$. Clearly the map $S$ is onto. Here is our
second main result.

\begin{theorem}\label{T:main-res2}
Let $n:=\dim V \geq 4$ and $2\leq k\leq n-2$. There exists a closed
proper subspace of $C^\infty(Gr_k(V))$ with the following property:
if $f\in C^{\infty}(G_k(V))$ is such that $\phi_f$ extends by
continuity to $\ck(V)$, then $S(f)$ belongs to this closed proper
subspace. This implies that the subset of $f\in C^\infty(G_k(V))$
such that $\phi_f$ extends by continuity to $\ck(V)$ is contained in
a closed proper linear subspace.
\end{theorem}

\begin{remark}
This theorem already implies Theorem \ref{T:main-res1} for $2\leq
k\leq n-2$. Actually, as it will be seen from the proof, the closed
proper subspace of $C^\infty(Gr_k(V))$ mentioned in the theorem is
equal to the image of the so called cosine transform on the
Grassmannian $Gr_k(V)$. For $2\leq k\leq n-2$ it is known to be a
proper subspace by a result of Goodey, Howard, and Reeder
\cite{goodey-howard-reeder}; see also \cite{alesker-bernstein} for a
more precise description of this subspace. The proof of the
proposition uses the Klain imbedding \cite{klain2} of even
valuations to functions on the Grassmannian and the result by J.
Bernstein and the author \cite{alesker-bernstein} stating that the
image of the Klain imbedding is equal to the range of the cosine
transform.
\end{remark}

The case $k=1$ is not covered by Theorem \ref{T:main-res2}; here the
situation is more subtle since the range of the cosine transform on
$C^\infty(Gr_1(V))$ is equal to the whole space. Here we will
(easily) see that in order to prove Theorem \ref{T:main-res1} it
suffices to do that for $\dim V=3$. In this case we have a much more
precise statement. In order to formulate it, it will be necessary to
rewrite the construction of $\phi_f$ in a more invariant way without
using explicitly the Euclidean metric on $V$. Then the group
$GL_n(\RR)$ will act on all spaces. It will be necessary to consider
$f$ not as a function on $\gk$ but as a section of a certain
$GL_n(\RR)$-equivariant line bundle over $\gk$.

Let us describe this. We will do it for $\dim V\geq 3$ and any
$1\leq k\leq n-2$. First in metric free terms we identify $G_k(V)$
with the set of pairs $(E,l)$ where $E\in Gr_k(V)$ is a
$k$-dimensional linear subspace, and $l$ is a line contained in the
annihilator $E^\perp \subset V^*$.

Let $p:\gk\to Gr_k(V)$ denote the projection $p(E,l)=E$. Let $L\to
Gr_k(V)$ be the line bundle whose fiber over $E\in Gr_k(V)$ is equal
to the space of Lebesgue measures on $E$. Note that a choice of a
Euclidean metric on $V$ defines a trivialization of $L$.

Next let us denote by $|\omega_{G_k/Gr_k}|$ the line bundle over
$\gk$ of relative densities of $\gk$ over $Gr_k(V)$. Let us recall
its definition. Let $q\colon X\to Y$  be a smooth map between smooth
manifolds. Assume that $q$ is a submersion, i.e. its differential at
every point is onto. (In our case $X=G_k(V), Y=Gr_k(V), q=p$.) The
fiber over a point $x\in X$ of the line bundle of relative densities
$|\omega_{X/Y}|$ is, by definition, equal to the (1-dimensional)
space of Lebesgue measures on the tangent space at $x$ to the fiber
$p^{-1}(p(x))$. If $q$ is proper, then we have a canonical map
$C^\infty(X,|\omega_{X/Y}|)\to C^\infty(Y)$ given by integration
along the fibers of $p$.

Note also that a choice of a Euclidean metric on $V$ defines a
trivialization of $|\omega_{G_k/Gr_k}|$ (indeed a Euclidean metric
on $V$ defines a Haar measure on each fiber of $p$). Let
$$M:=p^*L\otimes |\omega_{G_k/Gr_k}|.$$ Thus a choice of a Euclidean
metric on $V$ defines a trivialization of $M$, and hence an
isomorphism
\begin{eqnarray}\label{E:identification}
C^\infty(G_k(V),M)\simeq
C^\infty(G_k(V)).
\end{eqnarray}

In this notation and with the isomorphism (\ref{E:identification}),
the above construction of $\phi_f$ can be rewritten as follows. Let
$f\in C^\infty(\gk,M)$. Let $P\in \cp(V),\, F\in \cf_k(P)$. Then
$\int_{\gamma_F}f$ is a Lebesgue measure on $\bar F$. \footnote{To
be more precise, $\gamma_F$ is a subset of the manifold
$\PP_+(E^\perp)$ of oriented lines in $E^\perp$, and here we
integrate $f$ over the image of $\gamma_F$ in $\PP(E^\perp)$ under
the natural two sheeted covering $\PP_+(E^\perp)\to \PP(E^\perp)$
given by forgetting the orientation.} Hence
$\int_F\left(\int_{\gamma_F}f\right)$ is a number. Then
\begin{eqnarray}\label{E:phi-f}
\phi_f(P)= \sum_{F\in\cf_k(P)}\int_F\left(\int_{\gamma_F}f\right).
\end{eqnarray}
Set $$\Xi :=\{f\in C^\infty(\gk,M)|\, \phi_f \mbox{ extends by
continuity to } \ck(V)\}.$$ Thus under the identification
(\ref{E:identification}) the subspace $\Xi$ is isomorphic to the
space of smooth functions $f$ on $ G_k(V)$ such that $\phi_f$
extends by continuity to $\ck(V)$.

In Proposition \ref{1-2} we will show
\begin{proposition}\label{P:closedness}
$\Xi$ is a closed $GL(V)$-invariant linear subspace of
$C^\infty(\gk,M)$.
\end{proposition}

The next main result of this paper is the following one; later on it
will be formulated more explicitly.

\begin{theorem}\label{T:dim-3-main}
Let $n=3,\, k=1$. The natural representation of $GL_3(\RR)$ in
$C^\infty(G_1(V),M)$ has length three. The (closed
$GL_3(\RR)$-invariant) subspace $\Xi$ has length two. Moreover $\Xi$
coincides with the image of an explicit $GL_3(\RR)$-equivariant map
closely related to the construction of valuations using integration
with respect to the normal cycle (described below).
\end{theorem}

\begin{remark}
(1) As stated, Theorem \ref{T:dim-3-main} does not provide a unique
representation theoretical characterization of the subspace $\Xi$.
However with a little bit of extra work it can be shown that $\Xi$
is the only $GL_3(\RR)$-equivariant closed subspace of
$C^\infty(G_1(V),M)$ of length two whose irreducible subquotients
have the minimal (and equal to each other) Gelfand-Kirillov
dimension among the three irreducible subquotients of the whole
space (the third irreducible subquotient has strictly larger
Gelfand-Kirillov dimension). We will not pursue this charaterization
in this paper, but people familiar with the Gelfand-Kirillov
dimension can easily deduce the above characterization from our
arguments.

(2) For any $n$ and $k$ the recent preprint \cite{hinderer-hug-weil}
contains a sufficient condition for a section to belong to $\Xi$
(see Theorem 4.1 there).
\end{remark}

The proof of this theorem uses the Klain imbedding theorem,
construction of valuations using integration with respect to the
normal cycle, and a detailed study of the representation of
$GL_3(\RR)$ in $C^\infty(G_1(\RR^3),M)$; the latter step uses the
Beilinson-Bernstein localization theorem.

\hfill

The paper is organized as follows. In Section \ref{sec-valuations}
we prove the main results on valuations except the case $n=3,k=1$.
In Section \ref{S:case-3-1} we prove the case $n=3,k=1$ modulo some
representation theoretical computations which are postponed to
Section \ref{sec-gl3}. In Section \ref{sec-localization} we remind
some representation theoretical background; that section does not
contain new results. In Section \ref{sec-gl3} we make the
computations with representations of $GL_3(\RR)$; the main result of
this section is Theorem \ref{4-3}.

\hfill

{\bf Acknowledgements.} I am very grateful to S. Sahi who has
pointed out a mistake in an earlier version of the paper. I thank F.
Schuster for his numerous remarks on the first version of the paper.
I thank D. Gourevitch, D. Hug, and W. Weil for useful discussions.

\section{Proof of main results on
valuations when $n>3$.}\label{sec-valuations} We keep the notation
of the introduction.


\begin{proposition}\label{1-2}
The space $\Xi$ is a closed $GL(V)$-invariant linear subspace of
$C^\infty(\gk,M)$.
\end{proposition}

The linearity and $GL(V)$-invariance of $\Xi$ are obvious. In order
to prove that $\Xi$ us closed we will need some preparations.

\def\vek{Val_k^{ev}(V)}
\def\veks{Val_k^{ev,sm}(V)}
Let us denote by $Val_k^{ev}(V)$ the space of translation invariant
continuous even $k$-homogeneous valuations on $\ck(V)$. Equipped
with the topology of uniform convergence on compact subsets of
$\ck(V)$, $\vek$ is a Banach space. Let us recall the construction
of the Klain imbedding
$$K:\vek \hookrightarrow C(Gr_k,L)$$
which is a continuous linear map commuting with the natural action
of $GL(V)$.  Let $\phi\in \vek$. For any subspace $E\in Gr_k(V)$ the
restriction of $\phi|_E$ is a $k$-homogeneous valuation. By a result
of Hadwiger \cite{hadwiger-book}, $\phi|_E$ is a Lebesgue measure on
$E$. Let us denote by $K$ the map $\phi \mapsto[E\mapsto \phi|_E]$.
Thus $K:\vek \to C(Gr_k(V),L)$. A deep result due to D. Klain
\cite{klain2} (heavily based on \cite{klain1}) says that $K$ is
injective.

Let us denote by $\veks$ the subspace of {\itshape smooth}
$k$-homogeneous even valuations (a valuation $\phi\in Val(V)$ is
called smooth if the map $GL(V)\to Val(V)$ given by $g\mapsto
g(\phi)$ is $C^\infty$-differentiable). $\veks$ has a natural
structure of a Fr\'echet space. Moreover
$$K:\veks\hookrightarrow C^\infty(Gr_k(V),L).$$
The Casselman-Wallach theorem (see Theorem \ref{casselman-wallach}
below) implies the following claim.
\begin{claim}\label{1-3}
The Klain imbedding
$$K:\veks\hookrightarrow C^\infty(Gr_k(V),L)$$
has a closed image and induces an isomorphism of $\veks$ onto its
image as topological vector spaces.
\end{claim}

\def\ci{C^\infty}
Let us denote by $T$ the map $\Xi\to \vek$ given by $f\mapsto
\phi_f$. Note also that we have the canonical map
\begin{eqnarray}\label{E:map-S}
S:\ci(\gk,M)\to \ci(Gr_k(V),L)
\end{eqnarray}
given by the integration along the fibers of the projection
$p:\gk\to Gr_k(V)$. The map $S$ was defined in the introduction
under a choice of Euclidean metric. To rewrite it in invariant terms
let us recall that by definition $M=p^*L\otimes \omega_{G_k/Gr_k}$.
Then for any $f\in C^\infty(G_k(V),M)$ and any $E\in Gr_k(V)$, the
restriction $f|_{p^{-1}(E)}$ is a section over $p^{-1}(E)$ of the
line bundle $|\omega_{p^{-1}(E)}|\otimes L|_E$.  Hence one has
$\int_{p^{-1}(E)}f\in L|_E$. Now it is clear that $S$ is a
continuous linear map commuting with the natural action of $GL(V)$.

The next lemma is obvious.
\begin{lemma}\label{1-4}
The composition
$$K\circ T:\Xi\to C(Gr_k,L)$$
is equal to $S|_\Xi$.
\end{lemma}

\begin{lemma}\label{1-5}
For any $f\in \Xi$ the valuation $\phi_f$ is smooth, i.e.
$$T(\Xi)\subset \veks.$$
\end{lemma}
{\bf Proof.} By Lemma \ref{1-4} for any $f\in \Xi$, $K(\phi_f)\in
\ci(Gr_k(V),L)$. Let us show that this implies that $\phi_f$ is
smooth. By Claim \ref{1-3}, $K(\veks)$ is a closed subspace of
$\ci(Gr_k(V),L)$. Its closure $\overline{K(\veks)}$ in
$C(Gr_k(V),L)$ contains $K(\vek)$. It is easy to see that
$$\overline{K(\veks)}\cap \ci(Gr_k(V),L)= K(\veks).$$ Hence there
exists $\psi\in \veks$ such that $K(\psi)=K(\phi_f)$. Since $K$ is
injective on $\vek$ we get $\psi=\phi_f$. \qed

\hfill

{\bf Proof} of Proposition \ref{1-2}. It remains to show that
$\Xi$ is a closed subspace of $\ci(\gk,M)$. First let us show that
the map $T:\Xi\to \veks$ is continuous when $\Xi$ is equipped with
the topology induced from $\ci(\gk,M)$. By Claim \ref{1-3} it is
enough to show that $K\circ T$ is continuous. By Lemma \ref{1-4}
$K\circ T=S$ which is obviously continuous.

Now let us assume $f\in \overline{\Xi}$. Let $\{f_N\}\subset \Xi$ be
a sequence such that $f_N\to f$. The sequence $\{Tf_N\}\subset
\veks$ is a Cauchy sequence due to the continuity of $T$.  Hence
there exists $\phi\in \veks $ such that $Tf_N=\phi_{f_N} \to \phi$.
Hence $(K\circ T)(f_N)\to K(\phi)$. Also $Sf_N\to Sf$. Hence
$Sf=K(\phi)$. By construction, it is also clear that for any
polytope $P\in \cp(V)$ one has $\phi_{f_N}(P)\to \phi_f(P)$. But
also $\phi_{f_N}(P)\to \phi(P)$. Hence $\phi_f(P)=\phi(P)$ for any
$P\in \cp(V)$. This means that $\phi_f$ admits an extension by
continuity to $\ck(V)$. Hence $f\in \Xi$. \qed

\begin{theorem}\label{1-6}
$$T(\Xi)=\veks.$$
\end{theorem}
{\bf Proof.} Recall that $T:\Xi\to \veks$ is a continuous
$GL(V)$-equivariant map. Moreover it is non-zero by Example
\ref{example1}. By the Irreducibility Theorem \cite{alesker-gafa}
the space $\veks$ is $GL(V)$-irreducible. Hence $Im T$ is dense in
$\veks$. By Proposition \ref{1-2}, $\Xi$ is a closed subspace of
$\ci(\gk,M)$. Hence by the Casselman-Wallach theorem (see Theorem
\ref{casselman-wallach} below) $Im T$ is closed. Hence $Im T=\veks$.
\qed

\begin{theorem}\label{1-7}
The subspace $S(\Xi)\subset \ci(Gr_k(V),L)$ coincides with the
range of the cosine transform.
\end{theorem}
{\bf Proof.} Recall that by Lemma \ref{1-4} $S|_\Xi=K\circ T$. Thus
$S(\Xi)=(K\circ T)(\Xi)=K(\veks)$ where the last equality follows
from Theorem \ref{1-6}. But $K(\veks)$ coincides with the range of
the cosine transform by \cite{alesker-bernstein}. \qed

\hfill

Theorem \ref{1-7} immediately implies Theorem \ref{T:main-res2}.

So far we have proven Theorem \ref{T:main-res2}. It implies Theorem
\ref{T:main-res1} in all cases except when $k=1$ and $n\geq 3$. Let
us deduce the case $k=1,n\geq 3$ from the case $k=1,n=3$. Let $V$ be
an $n$-dimensional vector space, $n>3$. Let us fix a 3-dimensional
subspace $W\subset V$. Let $f\in C^\infty(G_1(V))$ be an arbitrary
smooth function. Let $\phi_f$ be the corresponding valuation on
polytopes $\cp(V)$. Let us describe explicitly the restriction of
$\phi_f$ to $\cp(W)\subset \cp(V)$ in terms intrinsic to $W$. We
will construct a function $g\in C^\infty(G_1(W))$ such that the
restriction of $\phi_f$ to $\cp(W)$ is equal to $\phi_g$.

The function $g$ can be described as follows. Let $p_{W^\perp}\colon
V\to W^\perp$ denote the orthogonal projection to $W^\perp$. Then
\begin{eqnarray}\label{E:f-g}
g(E,l)=\kappa_n\int_{m\in S(l\oplus
W^\perp)}f(E,m)\sqrt{1-|p_{W^\perp}(m)|^2}dm, \end{eqnarray} where
the integration is over the unit sphere of the space $l\oplus
W^\perp$, $dm$ is rotation invariant Haar measure on this sphere,
$|\cdot|$ under the square root is the Euclidean norm, and
$\kappa_n\ne 0$ is a normalizing constant depending on $n$ only. We
leave it to the reader to verify this elementary formula.

Clearly the map $f\mapsto g$ defined by (\ref{E:f-g}) is a
continuous linear operator $$U\colon C^\infty(G_1(V))\to
C^\infty(G_1(W)).$$ It is easy to see that this map is onto.

Let $g$ be an $O(3)$-finite function on $G_1(W)$ such that the
valuation $\phi_g$ does not extend by continuity to all convex
compact subsets of $W$; for the moment we just assume its existence.
It suffices to show that there exists an $O(n)$-finite function $f$
on $G_1(V)$ such that $U(f)=g$. Then we would have that the
restriction of $\phi_f$ to $\cp(W)$ is equal to $\phi_g$; this
implies that $\phi_f$ does not extend by continuity to all convex
compact subsets of $V$.

Notice that existence of such an infinitely smooth function $f$
follows from the surjectivity of $U$ on infinitely smooth functions.
Next we have the obvious imbedding $O(3)\subset O(n)$. The map $U$
commutes with the action of $O(3)$ on both spaces. Clearly any
$O(n)$-finite function is $O(3)$-finite. It follows that since
$O(n)$-finite functions are dense in $C^\infty(G_1(V))$, their image
under $T$ is dense in $C^\infty(G_1(W))$, and any function in this
image is $O(3)$-finite. But any dense linear $O(3)$-invariant
subspace of $C^\infty(G_1(W))$ consisting of $O(3)$-finite functions
is equal to the space of {\itshape all} $O(3)$-finite functions on
$C^\infty(G_1(V))$; this easily follows from the fact that the
natural representation of $O(3)$ in functions on $G_1(W)$ has finite
multiplicities due to transitivity of the action on $G_1(W)$. This
finishes the proof of Theorem \ref{T:main-res1}.

The following lemma will be needed in the next section.

\begin{lemma}\label{1-8}
Let $n\geq 3,\, 1\leq k\leq n-1$. Then
$$\Xi \cap Ker S\ne Ker S.$$
\end{lemma}
{\bf Proof.} Let us fix a Euclidean metric on $V$, and let us
identify $C^\infty(\gk, M)$ with $C^\infty(\gk)$. In order to prove
the lemma it is enough to construct a function $f\in Ker S$ and a
simplex $\Delta$ such that $\phi_f(\Delta)\ne 0$. Let us fix an
arbitrary $n$-simplex $\Delta$. Let $F_1,\dots, F_t$ be all of its
$k$-faces (then $t={n+1 \choose k+1}$). Let $A_i:=p^{-1}(\bar
F_i)\subset \gk$. Thus each $A_i$ is isomorphic to $\RR
\PP^{n-k-1}$. Let us choose an arbitrary function $g\in
C^\infty(\gk)$ such that $g|_{A_i}\equiv 0$ for $i>1$, $\int_{\gamma
_{\bar F_1}}g\ne 0$, and $\int_{A_1} g=0$. Then clearly
$\phi_g(\Delta)\ne 0$. Let $h:=S(g)\in C^\infty(Gr_k(V))$. Let
$f:=g-p^*g$. Then clearly $f\in Ker S$, $f|_{A_i}\equiv g|_{A_i}$
for any $i=1,2,\dots, t$, and hence
$\phi_f(\Delta)=\phi_g(\Delta)\ne 0$. Hence $f\notin \Xi$. This
proves the lemma. \qed

\section{The case $n=3,\, k=1$.}\label{S:case-3-1}
The main result of this section, Theorem \ref{T:main-3-1}, is stated
below. Let us start with some technical preparations.

Let us study the natural representation of the group $GL(V)$ in the
space $C^\infty(\gk,M)$. Let us denote by $P$ the parabolic subgroup
of $GL(V)$ with the Levi component $GL_k\times GL_{n-k-1}\times
GL_1$. Let $\chi:P\to \CC^*$ be the character of $P$ given by
$$\chi \left(\left[ \begin{array}{ccc}
                       A&*&*\\
                       0&B&*\\
                       0&0&C
                       \end{array}\right]\right)=
                       |\det A|^{-1}\cdot |\det B|\cdot |C|^{-(n-k-1)}$$
where $A\in GL_k,\, B\in GL_{n-k-1},\, C\in GL_{1}$.

The following proposition is straightforward.
\begin{proposition}\label{1-9}
The natural representation of $GL(V)$ in the space $C^\infty(\gk,M)$
is equal to the induced representation $Ind_P^G \chi$ (the induction
is not unitary).
\end{proposition}

Let us study in greater detail the case $n=3,\, k=1$. Then the
corresponding parabolic subgroup $P$ is the subgroup of upper
triangular matrices of size 3, and
\begin{eqnarray}\label{E:def-char-chi}
\chi \left(\left[ \begin{array}{ccc}
                       a&*&*\\
                       0&b&*\\
                       0&0&c
                       \end{array}\right]\right)=|a|^{-1}|b||c|^{-1}.
\end{eqnarray}

In Section \ref{sec-gl3} below we will prove the following
representation theoretical result.

\begin{proposition}\label{1-11}
Let $n=3,\, k=1$. The representation $Ind_P^G\chi$ is of length 3 as
a $GL(V)$-representation. Moreover the irreducible subquotients are
pairwise non-isomorphic.
\end{proposition}

\hfill

Recall from Section \ref{sec-valuations} that we have the map
(integration along the fibers)
$$S\colon C^\infty(G_1,M)\to C^\infty(Gr_1(V),L).$$
We will need a lemma.
\begin{lemma}\label{L:length-ker-S}
$Ker S$ has length 2 as $GL(V)$-module.
\end{lemma}
{\bf Proof.} By Theorem \ref{1-7} the subspace $S(\Xi)\subset
C^\infty(Gr_1(V),L)$ coincides with the range of the cosine
transform. But in this case the cosine transform is an isomorphism
(see e.g. \cite{semyanistyi}). Moreover the source and the target
spaces of the cosine transform are irreducible
$GL(V)$-representations (this is well known; see e.g.
\cite{howe-lee}). Hence $Im(S)=C^\infty(Gr_1(V),L)$ is irreducible.
This and Proposition \ref{1-11} imply that the length of $Ker(S)$ is
equal to 2. \qed

\hfill

Let us remind an important construction of continuous valuations
called integration with respect to the normal cycle; it will be
necessary for the formulation of the main result of this section. For a
convex compact subset $K\subset V$ and a point $x\in K$ let us
define a tangent cone $T_xK$ to $K$ at $x$:
\begin{eqnarray*}
T_xK:=\overline{\cup_{t>0}t(K-x)}.
\end{eqnarray*}
Clearly $T_xK$ is a convex cone in $V$. Define the dual cone
$$(T_xK)^o:=\{\xi\in V^*|\, <\xi,v>\leq 0\mbox{ for all } v\in
T_xK\}.$$ Let us denote by $\PP_+(V^*)$ the oriented
projectivization of $V^*$, namely the set of oriented one
dimensional linear subspaces of $V^*$. Another way to view it is as
the quotient space $(V^*\backslash \{0\})/\RR_{>0}$, where
$\RR_{>0}$ acts of $V^*$ by usual multiplication.

Now let us define the normal cycle $N(K)$ as a subset of
$\PP_+(V^*)\times V$:
$$N(K):=\cup_{x\in K}(T_xK)^o/\RR_{>0}.$$
It is well known that $N(K)$ is a compact set of Hausdorff dimension
$n-1$ ($n=\dim V$). It is easy to see that an orientation of $V$
induces an orientation of $N(K)$.

\hfill

To construct a continuous valuation, let us fix
$$\ome\in
C^\infty(\PP_+(V^*),\Ome^{n-1})\otimes or(V),$$ where $\Omega^{n-1}$
denotes the vector bundle whose sections are $n-1$ forms, and
$or(V)$ is the orientation line of $V$. (Remind the definition of $or(V)$: 
this is the one dimensional space of $\CC$-valued functions $H$ on the set of all bases of $V$
satisfying the property $H(gB)=sgn(g)H(B)$ for any $g\in GL(V)$ and any basis $B$). 
Then define a functional on
$\ck(V)$:
$$K\mapsto \int_{N(K)}\ome.$$
It is well known (see e.g \cite{alesker-fu}, Section 2.1) that this
functional is a continuous valuation on $\ck(V)$. It is clear that
if $\ome$ is translation invariant with respect to $V$ then the
corresponding valuation is translation invariant. The space of
translation invariant $\ome$'s can be written
$$\oplus_{k=0}^{n-1}
C^\infty(\PP_+(V^*),\Ome^{n-k-1})\otimes\wedge^k V^*\otimes or(V).$$ Elements of
$C^\infty(\Omega^{n-1-k}(\PP_+(V^*))\otimes \wedge^kV^*\otimes
or(V))$ define $k$-homogeneous valuations. Even valuations are
determined by even sections of the latter space, when the parity on
them is induced by the involution on $\PP_+(V^*)$ of the change of
orientation. This space of even sections over $\PP_+(V^*)$ can be
viewed as the space of all sections of an appropriate vector bundle
over the projective space $\PP(V^*)$ of non-oriented lines in $V^*$.
It is not hard to see that this vector bundle is
$$C^\infty(\PP(V^*),\Ome^{n-k-1}\otimes\eps_{n-k})\otimes \wedge^kV^*\otimes or(V),$$
where $\eps_{n-k}$ is the trivial line bundle for
$n-k$ even, and for $n-k$ odd its fiber over $l\in \PP(V^*)$ is
equal to the orientation line $or(l)$ of $l$.

\hfill

Let us introduce more notation. Let $\cp Val^+_k$ denote the space
of all valuations on convex compact polytopes which are
$k$-homogeneous and even (without any continuity assumptions). Let
$$\cf\colon C^\infty(G_k(V),M)\to \cp Val_k^+$$
be given by $f\mapsto \phi_f$ as in (\ref{E:phi-f}). Let
$$\cg\colon C^\infty(\PP(V^*), \Ome^{n-k-1}\otimes \eps_{n-k})
\otimes\wedge^kV^*\otimes or(V) \to \cp Val^+_k$$ be the map
$\omega\mapsto [K\mapsto \int_{N(K)}\ome]$. The image of $\cg$ is
contained in continuous valuations $Val^+_k$. Clearly $\cf,\cg$ are
$GL(V)$-equivariant.
\begin{proposition}
There exists a $GL(V)$-equivariant map $$\ch\colon
C^\infty(\PP(V^*), \Ome^{n-k-1}\otimes \eps_{n-k})
\otimes\wedge^kV^*\otimes or(V)\to C^\infty(G_k(V),M)$$ such that
$\cg=\cf\circ \ch$.
\end{proposition}

{\bf Proof.} Let us construct $\ch$ explicitly. Fix $$\ome\in
C^\infty(\PP(V^*), \Ome^{n-k-1}\otimes \eps_{n-k})
\otimes\wedge^kV^*\otimes or(V).$$ Let $E\in Gr_k(V)$, and let
$l\subset E^\perp$ be a line. Restriction of $\ome$ to
$\PP(E^\perp)\subset \PP(V^*)$ is en element of
$C^\infty(\PP(E^\perp), \Ome^{n-k-1}\otimes \eps_{n-k})
\otimes\wedge^kV^*\otimes or(V)$. The canonical map $V^*\to E^*$
defines a linear map $\wedge^kV^*\to \wedge^k E^*$. Hence $\ome$
defines an element in
$$C^\infty(\PP(E^\perp), \Ome^{n-k-1}\otimes \eps_{n-k})
\otimes\wedge^k E^*\otimes or(V).$$ Taking the value of this
section over $l\in \PP(E^\perp)$ we get an element in
\begin{eqnarray}\label{E:space}
\wedge^{n-1-k}T_l^*\PP(E^\perp)\otimes \wedge^k E^*\otimes
or(V)\otimes \eps_{n-k}|_l.
\end{eqnarray}
Let us show that the last space is canonically isomorphic to
$$|\ome_{G_k/Gr_k}|\big|_l\otimes L|_E.$$ Clearly
\begin{eqnarray}\label{E:space2}
|\ome_{G_k/Gr_k}|\big|_l=\wedge^{n-1-k}T^*_l\PP(E^\perp)\otimes
or(T^*_l\PP(E^\perp)),\\\label{E:space3} L|_E=\wedge^k E^*\otimes
or(E).
\end{eqnarray}
Comparing (\ref{E:space2})-(\ref{E:space3}) with (\ref{E:space}) we
see that we have to construct a canonical isomorphism
\begin{eqnarray}\label{E:canon-iso}
or(T^*_l(\PP(E^\perp))\otimes or(E)=or(V)\otimes \eps_{n-k}|_l.
\end{eqnarray}
To prove this, we will use the following general canonical
isomorphisms
\begin{eqnarray*}
or(X^*)=or(X),\\
or(X/Y)=or(X)\otimes or(Y),\\
or(X\otimes Y)=or(X)^{\otimes\dim Y}\otimes or(Y)^{\otimes \dim X}.
\end{eqnarray*}

Then the left hand side in (\ref{E:canon-iso}) is equal to
\begin{eqnarray*}
or(Hom(l,E^\perp/l)^*)\otimes or(E)=\\
or(l\otimes (E^\perp/l)^*)\otimes or(E)=\\
or(l)^{\otimes(n-k-1)}\otimes or(E^\perp)^*\otimes or(l)\otimes or(E)=\\
or(l)^{\otimes(n-k)}\otimes or(V/E)\otimes
or(E)=or(l)^{\otimes(n-k)}\otimes or(V)=\\
or(V)\otimes \eps_{n-k}|_l.
\end{eqnarray*}

Thus the isomorphism (\ref{E:canon-iso}) is constructed. Hence the
map $\ch$ is constructed. It is easy to see that it satisfies
$$\cg=\cf\circ \ch.$$
\qed

\hfill

In the rest of this section we will assume that $n=3,\, k=1$. The
following theorem is the main result of this section.
\begin{theorem}\label{T:main-3-1}
Let $n=3,\, k=1$.

(1) The subspace $\Xi\subset C^\infty(G_1(V),M)$ is equal to the
image of $\ch$.

(2) The length of $\Xi$ as a $GL(V)$-module is equal to 2, while the
length of the whole space $C^\infty(G_1(V),M)$ is equal to 3, and
all irreducible subquotients are pairwise non-isomorphic.
\end{theorem}

\begin{remark}
Part (1) of Theorem \ref{T:main-3-1} provides some kind of analytic
(or geometric) description of $\Xi$. Part (2) gives some
complimentary information of the subspace $\Xi$: for example it may
take some extra work to show that part (1) implies that $\Xi$ does
not coincide with the whole space. In fact, the proof of part (1)
uses part (2).

Part (2), as it is stated, does not provide a unique
characterization of the subspace $\Xi$ inside $C^\infty(G_1(V),M))$.
However with a little bit of extra work one can obtain the following
slightly more precise statement: $\Xi$ is the only closed
$GL(V)$-invariant subspace of length 2 whose irreducible
subquotients have minimal Gelfand-Kirillov dimensions among
irreducible subquotients of $C^\infty(G_1(V),M)$. (In fact, one can
show that the Gelfand-Kirillov dimensions of the three irreducible
subquotients of the whole space are equal to 2,2, and 3.)
\end{remark}

{\bf Proof of Theorem \ref{T:main-3-1}.} First it will be useful to
compute the kernel of $\ch$. Recall that for $n=3,k=1$ we have
$$\ch\colon C^\infty(\PP(V^*),\Ome^1)\otimes V^*\otimes or(V)\to
C^\infty(G_1(V),M).$$ Let us fix $\ome\in
C^\infty(\PP(V^*),\Ome^1)\otimes V^*\otimes or(V)$. Clearly $\ome\in
Ker(\ch)$ if and only if for any 1-dimensional linear subspaces
$E\subset V$ and $l\subset E^\perp (\subset V^*)$ the image of
$\ome$ vanishes under the map
\begin{eqnarray}\label{E:map1}
T^*_l(\PP(V^*))\otimes V^*\otimes or(V)\to
T^*_l(\PP(E^\perp))\otimes E^*\otimes or(V)
\end{eqnarray}
which is just the tensor product of the natural restriction maps
\begin{eqnarray*}
T^*_l(\PP(V^*))\to T^*_l(\PP(E^\perp)),\\
V^*\to E^*,
\end{eqnarray*}
and of the identity map on $or(V)$.

Let us fix $l\subset V^*$ and describe explicitly the intersection
of kernels of the maps (\ref{E:map1}) over all lines $E\subset
l^\perp$. Thus we have
\begin{eqnarray*}
T^*_l(\PP(V^*))=l\otimes(V^*/l)^*=l\otimes l^\perp,\\
T^*_l(\PP(E^\perp))=l\otimes(E^\perp/l)^*=l\otimes (l^\perp/E).
\end{eqnarray*}

Then (\ref{E:map1}) becomes a map
\begin{eqnarray}\label{E:map2}
l\otimes (l^\perp\otimes V^*)\otimes or V\to
l\otimes((l^\perp/E)\otimes E^*)\otimes or(V),
\end{eqnarray}
where the map is tensor product of the identity map on the first
copy of $l$ and on $or(V)$, of the canonical quotient map
$l^\perp\to l^\perp/E$, and of the obvious natural map $V^*\to E^*$.

It suffices to describe the intersection of the kernels of maps
\begin{eqnarray}\label{E:map3}
l^\perp\otimes V^*\to (l^\perp/E)\otimes E^*,
\end{eqnarray}
when $l$ is fixed, and $E$ runs over all lines $E\subset l^\perp$.
The map (\ref{E:map3}) can be identified with the map
\begin{eqnarray}\label{E:map4}
Hom(V,l^\perp)\to Hom(E,l^\perp/E)
\end{eqnarray}
given by the composition $\phi\mapsto p_E\circ\phi\circ i_E$, where
$i_E\colon E\to V$ is the identity imbedding, and $p_E\colon
l^\perp\to l^\perp/E$ is the quotient map. The kernel of
(\ref{E:map4}) consists of maps $\phi\colon V\to l^\perp$ such that
$\phi(E)\subset E$. Hence the intersection of kernels of
(\ref{E:map4}) over all $E\subset l^\perp$, when $L$ is fixed,
consists of $\phi\colon V\to l^\perp$ such that any line $E\subset
l^\perp$ is mapped to $E$, or equivalently the restriction of $\phi$
to $l^\perp$ is proportional to the identity map. Let us denote by
$X_l$ this space of maps tensored by $l\otimes or(V)$. Then $X_l$ is
the intersection of the kernels of (\ref{E:map2}), which is the same
as the intersection of kernels of (\ref{E:map1}) over all $E\subset
l^\perp$ with fixed $l$.

Now let us describe explicitly the quotient of
$T^*_l(\PP(V^*))\otimes V^*\otimes or(V)$ by $X_l$. Clearly
$$X_l\supset l\otimes Hom(V/l^\perp,l^\perp)\otimes or(V)=l\otimes(l^\perp\otimes
l)\otimes or(V).$$ Then we have furthermore
\begin{eqnarray*}
(T^*_l(\PP(V^*))\otimes V^*\otimes or(V))/(l\otimes(l^\perp\otimes
l)\otimes or(V))=\\
(l\otimes(l^\perp\otimes V^*)\otimes or(V))/(l\otimes(l^\perp\otimes
l)\otimes or(V))=\\
l\otimes(l^\perp\otimes V^*/l)\otimes or(V)=\\
l\otimes End(V^*/l)\otimes or(V).
\end{eqnarray*}
We have to quotient out the last space by
$X_l/(l\otimes(l^\perp\otimes l)\otimes or(V))$. The result is
\begin{eqnarray}\label{E:map5}
l\otimes End_0(V^*/l)\otimes or(V),
\end{eqnarray}
where $End_0$ denotes the trace zero endomorphisms.

Let us denote by $N$ the bundle over $\PP(V^*)$ whose fiber over $l$
is equal to (\ref{E:map5}). Clearly $N$ is $GL(V)$-equivariant. Thus
so far we obtained that the map $\ch$ uniquely factorizes via a
linear map
\begin{eqnarray}\label{E:map6}
\bar\ch\colon C^\infty(\PP(V^*),N)\to C^\infty(G_1(V),M)
\end{eqnarray}
which is an injective $GL(V)$ equivariant map. As we have mentioned
previously, $Im(\ch)\subset \Xi$. Hence
\begin{eqnarray}\label{E:inclusion}
Im (\bar \ch)\subset \Xi.
\end{eqnarray}

It is easy to see that the composition
\begin{eqnarray*}
S\circ \bar\ch\colon C^\infty(\PP(V^*),N)\to C^\infty(Gr_1(V),L)
\end{eqnarray*}
is a non-zero $GL(V)$-equivariant map. Since the target is an
irreducible representation, the image of $S\circ \bar\ch$ is
everywhere dense. Furthermore the Casselmann-Wallach theorem (see
Theorem \ref{casselman-wallach} below) implies immediately that
$S\circ \bar \ch$ is onto, but we will not need this fact.


Let us show that
\begin{eqnarray}\label{E:kernel}
Im (\bar\ch)\cap Ker S\ne \{0\}.
\end{eqnarray}
A short proof of this claim was personally communicated to us by S.
Sahi. We reproduce here his argument.

Let us assume on the contrary that $Im \bar \ch\cap Ker S=\{0\}$.
Then
$$S\circ \bar \ch\colon C^\infty(\PP(V^*),N)\to
C^\infty(Gr_1(V),L)$$ is an isomorphism of the underlying
Harish-Chandra modules (and in fact of the actual spaces by the
Casselmann-Wallach theorem, though we will not need this fact here).
We will get a contradiction by comparing the $SO(3)$-actions of
them.

Let us fix a Euclidean metric and an orientation on $V$. Then the
space $C^\infty(Gr_1(V),L)$ can be $SO(3)$-equivariantly identified
with the space of smooth functions on $\RR\PP^2$. It is well known
that this space is a direct sum of irreducible representations of
$SO(3)$ with highest weights $2m$, $m\in \ZZ_{\geq 0}$, each
entering with multiplicity 1 exactly.

Now let us show that $C^\infty(\PP(V^*),N)$ contains some extra
irreducible $SO(3)$-representations. In the following computation
the symbol "$0$" on the bottom denotes the traceless part. We have
\begin{eqnarray*}
N|_l=l\otimes End_0(V^*/l)\otimes or(V)=\\
l\otimes[(V^*/l)^*\otimes (V^*/l)]_0\otimes
or(V)\overset{SO(3)}{\simeq}\\
l\otimes[(V^*/L)^{\otimes 2}]_0=\\
l\otimes (\wedge^2(V^*/l)\oplus [Sym^2(V^*/l)]_0).
\end{eqnarray*}
But $l\otimes\wedge^2(V^*/l)=\wedge^3
V^*\overset{SO(3)}{\simeq}\RR$. Thus we see that the $SO(3)$-module
$C^\infty(\PP(V^*),N)$ is a direct sum of the space of functions on
$\PP(V^*)$ (which coincides with the previous space) and another
infinite dimensional space. This implies  (\ref{E:kernel}).


Since $S\circ\bar \ch\ne 0$ and $Im (\bar \ch)\cap Ker S\ne 0$, it
follows that the length of $C^\infty(\PP(V^*),N)$ as $GL(V)$-module
is at least 2. The injectivity of $\bar \ch$ and (\ref{E:inclusion})
imply
\begin{eqnarray}\label{E:lengths}
2\leq \mbox{length}(C^\infty(\PP(V^*),N))\leq \mbox{length}(\Xi).
\end{eqnarray}

But by Proposition \ref{1-11} and Lemma \ref{1-8} we obtain
$$\mbox{length}(\Xi)<3.$$
Consequently
$$\bar\ch(C^\infty(\PP(V^*),N))=\Xi,$$
and both have length 2.

This implies Theorem \ref{T:main-3-1}.


\section{Reminder on representation theory.}\label{sec-localization} We recall the
Beilinson-Bernstein theorem on localization of $\mathfrak{g}
$-modules following \cite{beilinson-bernstein}, and the
Casselman-Wallach theorem \cite{casselman}, \cite{wallach}.

Let $G$ be a complex reductive algebraic group. Let $T$ denote a
Cartan subgroup of $G$. Let $B$ be a Borel subgroup of $G$
containing $T$. Let $\mathfrak{g}$ denote the Lie algebra of $G$,
$\mathfrak{t}$ the Lie algebra of $T$, and $\mathfrak{b}$ the Lie
algebra of $B$. Let $\mathfrak{n}$ denote the nilpotent radical of
$\mathfrak{b}$.

\begin{remark}
In our applications we will take $G=GL_3(\CC)$, $T$ will be the
subgroup of diagonal invertible matrices, $B$ will be the subgroup
of upper triangular invertible matrices, $\mathfrak{n}$ is the Lie
algebra of upper triangular matrices with zeros on the diagonal.
\end{remark}

Let $R(\mathfrak{t})\subset \mathfrak{t}^*$ be the set of roots of
$\mathfrak{t}$ in $\mathfrak{g}$. The set $R(\mathfrak{t})$ is
naturally divided into the set of roots whose root spaces are
contained in $\mathfrak{n}$ and its complement. The latter set is
denoted by $R^+(\mathfrak{t})$; it is equal to the set of roots of
$\mathfrak{t}$ in $\mathfrak{g}/\mathfrak{b}$. If $\alpha$ is a root
of $\mathfrak{t}$ in $\mathfrak{g}$ then the dimension of the
corresponding root subspace $\mathfrak{g}_{\alpha}$ is called the
multiplicity of $\alpha$. Let $\rho_{\mathfrak{b}}$ be the half sum
of the roots contained in $R^+(\mathfrak{t})$ counted with their
multiplicities.

For a root $\alpha$ let $\alpha^V\in \mathfrak{t}$ denote the
corresponding coroot (see e.g. \cite{goodman-wallach}, \S 2.4.2). We
say that $\lambda \in \mathfrak{t}^*$ is {\itshape dominant} if for
any root $\alpha \in R^+(t)$ we have $<\lambda, \alpha ^{V}>\ne -1,
-2, \dots$. We shall say that $\lambda \in t^*$ is {\itshape
regular} if for any root $\alpha \in R^+(t)$ we have $<\lambda,
\alpha ^{V}>\ne 0$.

For the definitions and basic properties of the sheaves of twisted
differential operators we refer to \cite{bien}. Here we will present
only the explicit description of the sheaf $\dl$ in order to agree
about the normalization.

\def\ug{U(\mathfrak{g})}
\def\ox{{\cal O}_X}
\def\tx{{\cal T}_X}
Let $X$ be the full flag variety of $G$ (then $X=G/B$). $X$ can be
identified with the variety of all Borel subalgebras of
$\mathfrak{g}$. Let $\ox$ denote the sheaf of regular functions on
$X$. Let $\ug$ denote the universal enveloping algebra of
$\mathfrak{g}$. Let $U^o$ be the sheaf $\ug \otimes _{\CC} \ox$, and
$\mathfrak{g}^o:=\mathfrak{g}\otimes _{\CC}\ox$. Let $\tx$ be the
tangent sheaf of $X$. We have the canonical morphism
$\alpha:\mathfrak{g}^o\to \tx$. Let also $\mathfrak{b}^o:=Ker
\,\alpha =\{\xi\in\mathfrak{g}^o|\,\xi_x\in \mathfrak{b}_x \forall
x\in X\},$ where $\mathfrak{b}_x\subset \mathfrak{g}$ denotes the
Borel subalgebra corresponding to $x\in X$. Let $\lambda\colon
\mathfrak{b}\to \CC$ be a linear functional which is trivial on
$\mathfrak{n}$ (thus $\lambda \in \mathfrak{t}^*$). Since
$\mathfrak{b}_x/[\mathfrak{b}_x,\mathfrak{b}_x]$ are canonically
isomorphic for different $x\in X$, $\lambda$ defines a morphism
$\lambda ^o\colon \mathfrak{b}^o\to \ox$. We will denote by ${\cal
D}_\lambda$ the sheaf of twisted differential operators
corresponding to $\lambda -\rho _{\mathfrak{b}}$, i.e. $\dl$ is
isomorphic to $U^o/{\cal I}_{\lambda}$, where ${\cal I}_{\lambda}$
is the two sided ideal generated by the elements of the form
$$\xi-(\lambda -\rho_{\mathfrak{b}})^o(\xi),$$ where $\xi$ is a
local section of $\mathfrak{b}^o$. Let $D_\lambda:=\Gamma (X,\dl)$
denote the ring of global sections of $\dl$. We have a canonical
morphism $U(\mathfrak{g})\to D_\lambda$. For the full flag variety
$X$ this map is onto (\cite{beilinson-bernstein}). Let us also
denote by ${\cal D}_{\lambda}-mod$ (resp. $D_\lambda -mod$) the
category of $\dl$- (resp. $D_{\lambda}-$) modules.

In this notation one has the following result due to Beilinson and
Bernstein \cite{beilinson-bernstein}.

\begin{theorem}[Beilinson-Bernstein]\label{2-1}
(1) If $\lambda \in \mathfrak{t}^*$ is dominant then the functor of
global sections $\Gamma\colon {\cal D}_{\lambda}-mod \to
U(\mathfrak{g})-mod$ is exact, namely $\Gamma$ maps exact sequences
to exact ones.

(2) If $\lambda \in \mathfrak{t}^*$ is dominant and regular then the
functor $\Gamma$ is also fully faithful, namely $\Gamma$ induces a
bijection on the set of morphisms between any two objects of ${\cal
D}_{\lambda}-mod$.
\end{theorem}
\begin{remark}
Let us notice that when $\lam=\rho_{\mathfrak{b}}$ we get the usual
(untwisted) ring $D$ and sheaf $\cd$ of differential operators. This
$\lam=\rho_{\mathfrak{b}}$ is regular and dominant, and hence the
Beilinson-Bernstein theorem is applicable.
\end{remark}

Note also that the functor $\Gamma$ can be considered as taking
values in the category $D_\lambda-mod$, rather than
$U(\mathfrak{g})-mod$. Considered in this way, $\Gamma$ has a left
adjoint functor, called the Beilinson-Bernstein localization
functor, $\Delta: D_\lambda -mod\to \dl -mod$. It is defined as
$\Delta (M)=\dl \otimes_{D_\lambda} M$.

Part (1) of the next lemma was proved in \cite{beilinson-bernstein},
parts (2), (3) can be found in \cite{bien}, Proposition I.6.6.

\begin{lemma}\label{2-2}
Suppose $\Gamma: {\cal D}_{\lambda}-mod \to D_{\lambda}-mod$ is
exact. Then

(1) the localization functor $\Delta:D_{\lambda}-mod\to\dl -mod$ is
the right inverse of $\Gamma$:
$$\Gamma\circ \Delta =Id;$$

(2) $\Gamma$ sends simple objects to simple ones or to zero.

(3) $\Gamma$ sends distinct simple objects to distinct ones or to
zero.
\end{lemma}

\hfill

Now let us discuss the Casselman-Wallach theorem.

\begin{definition}\label{part1-rep-1} Let $\pi$ be a continuous representation of a Lie group $G_0$ in a
Fr\'echet space $F$. A vector $\xi \in F$ is called $G_0$-smooth if
the map $g\mapsto \pi(g)\xi$ is an infinitely differentiable map
from $G_0$ to $F$.
\end{definition}
It is well known (see e.g. \cite{wallach}, Section 1.6) that the
subset $F^{sm}$ of smooth vectors is a $G_0$-invariant linear
subspace dense in $F$. Moreover it has a natural topology of a
Fr\'echet space (which is stronger than the topology induced from
$F$), and the representation of $G_0$ in $F^{sm}$ is continuous.
Moreover all vectors in $F^{sm}$ are $G_0$-smooth.

Let $G_0$ be a real reductive group. Assume that $G_0$ can be
imbedded into the group $GL_N(\RR)$ for some $N$ as a closed
subgroup invariant under the transposition. Let us fix such an
imbedding $p:G_0\hookrightarrow GL_N(\RR)$. (In our applications
$G_0$ will be either $GL_n(\RR)$ or a direct product of several
copies of $GL_n(\RR)$.) Let us introduce a norm $|\cdot |$ on $G_0$
as follows:

$$|g|:=\max\{||p(g)||,||p(g^{-1})||\}$$
where $||\cdot||$ denotes the usual operator norm in $\RR^N$.
\begin{definition}
Let $\pi$ be a smooth representation of $G_0$ in a Fr\'echet space
$F$ (namely $F^{sm}=F$). One says that this representation has
{\itshape moderate growth} if for each continuous semi-norm
$\lambda$ on $F$ there exists a continuous semi-norm $\nu_\lambda$
on $F$ and $d_{\lambda}\in \RR$ such that
$$\lambda(\pi(g)v)\leq |g|^{d_\lambda}\nu_{\lambda}(v)$$
for all $g\in G,\, v\in F$.
\end{definition}

The proof of the next lemma can be found in \cite{wallach}, Lemmas
11.5.1 and 11.5.2.
\begin{lemma}\label{part1-wallach}
(i) If $(\pi,G_0,H)$ is a continuous representation of $G_0$ in a
Banach space $H$, then $(\pi,G_0,H^{sm})$ has moderate growth.

(ii) Let $(\pi, G_0,V)$ be a representation of moderate growth. Let
$W$ be a closed $G_0$-invariant subspace of $V$. Then the
representations of $G_0$ in $W$ and $V/W$ have moderate growth.
\end{lemma}

Recall that a continuous Fr\'echet representation $(\pi,G_0,F)$ is
said to have {\itshape finite length} if there exists a finite
filtration
$$0=F_0\subset F_1\subset \dots\subset F_m=F$$
by $G_0$-invariant closed subspaces such that $F_i/F_{i-1}$ is
irreducible for any $i$, i.e. it does not have proper closed
$G_0$-invariant subspaces. The sequence of all consecutive quotients
$$F_1,F_2/F_1,\dots,F_m/F_{m-1}$$ is called the Jordan-H\"older
series of the representation $\pi$. It is well known (and easy to
see) that the Jordan-H\"older series of a finite length
representation is unique up to a permutation.

\begin{definition}
A Fr\'echet representation $(\rho,G_0,F)$ of a real reductive group
$G_0$ is called {\itshape admissible} if its restriction to a
maximal compact subgroup $K$ of $G_0$ contains an isomorphism class
of any irreducible representation of $K$ with at most finite
multiplicity. (Recall that a maximal compact subgroup of $GL_n(\RR)$
is the orthogonal group $O(n)$.)
\end{definition}

\begin{theorem}[Casselman-Wallach \cite{casselman}, \cite{wallach}]\label{casselman-wallach}
Let $G_0$ be a real reductive group. Let $(\rho,G_0,F_1)$ and
$(\pi,G_0,F_2)$ be smooth representations of moderate growth in
Fr\'echet spaces $F_1, F_2$. Assume in addition that $F_2$ is
admissible of finite length. Then any continuous morphism of
$G_0$-modules $f:F_1\to F_2$ has closed image.
\end{theorem}

\section{Representations of $GL_3(\RR)$}\label{sec-gl3}
\def\lch{{\cal L}^{\chi}}
Let $P$ be the Borel subgroup $GL_3(\RR)$ consisting of upper
triangular $3\times 3$ matrices. The positive roots are $(-1,1,0),
(-1,0,1),(0,-1,1)$. Hence the half sum of the positive roots is
$\rho =(-1,0,1)$.

The main goal of this section is to prove Proposition \ref{1-11}. We
will immediately see that Proposition \ref{1-11} is equivalent to
the following result.
\begin{theorem}\label{4-3}
Let $\psi(x,y,z)=|x|^{-1}$. Then the representation $Ind_P^G\psi$
(the induction is not normalized!) has length three, and all the
irreducible subquotients are pairwise non-isomorphic.
\end{theorem}
The equivalence of Proposition \ref{1-11} and Theorem \ref{4-3}
follows from an easy observation that the two induced
representations considered in these statements are associate, and
the well known fact (see e.g. \cite{vogan}, Proposition 4.1.20) that
associate representations have the same Jordan-H\"older series.

To prove Theorem \ref{4-3}, let us observe first that the character
$\psi$ is dominant and regular. Let us denote by $\cd_\psi$ the
corresponding sheaf of twisted differential operators on the flag
manifold $X$.

Now let us construct an $O(3,\CC)$-equivariant $\cd_\psi$-module
$\cm_\psi$ on $X$ such that the $U(\mathfrak{g})$-module of global
sections $\Gamma(X,\cm_\psi)$ coincides with the Harish-Chandra
module of $Ind_P^G\psi$. Let $U$ be the open $O(3,\CC)$-orbit in
$X$. Let $j\colon U\hookrightarrow X$ be the identity imbedding
(which is an affine morphism). Fix a flag $(E\subset F)\in U$. The
category of $O(3,\CC)$-equivariant coherent $\cd_\psi$-modules on
$U$ is equivalent to the category of finite dimensional
representations of the group of connected components of the
stabilizer of $(E,F)$ in $O(3,\CC)$. This stabilizer is equal to
$$O(1)\times O(1)\times O(1)\simeq \{\pm 1\}\times \{\pm 1\}\times \{\pm
1\}.$$ Let us consider the representation of this group given by
$$(A,B,C)\mapsto A.$$
Let $\cm_{\psi 0}$ be the corresponding $O(3,\CC)$-equivariant
$\cd_\psi$-module on $U$. (It is easy to see that $\cm_{\psi 0}$ is
a locally free rank one $\co_U$-module.) Define
$$\cm_\psi:=j_*\cm_{\psi 0}.$$
It is not hard to see that $\cm_\psi$ is the required
$\cd_\psi$-module.

By the Beilinson-Bernstein Theorem \ref{2-1} and Lemma \ref{2-2},
Theorem \ref{4-3} is equivalent to saying that $\cm_\psi$ has length
three with non-isomorphic irreducible subquotients in the category
of $O(3,\CC)$-equivariant $\cd_\psi$-modules.

The sheaf $\cd_\psi$ is the sheaf of differential operators acting
on algebraic sections of the $GL_3(\CC)$-equivariant line bundle
$\cl$ over $X$ which is defined as follows: for a complex flag
$(E\subset F)$ in $\CC^3$, the fiber of $\cl$ over $(E\subset F)$ is
equal to $E^*$. Then we have an equivalence of categories of
$O(3,\CC)$-equivariant $\cd_\psi$-modules and of
$O(3,\CC)$-equivariant $\cd$-modules (here $\cd$ denotes the sheaf
of usual untwisted differential operators on $X$) which is given by
$\cf\mapsto \cf\otimes_{\co_X}\cl^{-1}$, where $\co_X$ denotes the
sheaf of regular functions on $X$.

Thus we have to prove the following proposition.

\begin{proposition}\label{P:prop-legth2}
The $O(3,\CC)$-equivariant $\cd$-module
\begin{eqnarray}\label{M-def}
\cm:=\cm_\psi\otimes_{\co_X}\cl^{-1} \end{eqnarray} has length three
and all the irreducible subquotients are pairwise non-isomorphic.
\end{proposition}

Now let us describe explicitly the $O(3,\CC)$-orbits on the variety
of complete flags in $\CC^3$. Let us denote by $X$ the variety of
complete flags in $\CC^3$, namely the set of pairs $(E,F)$ where $E$
is a complex line, $F$ is a complex 2-plane, and $E\subset F$. The
following proposition is essentially well known and its proof is
left to the reader as an (easy) exercise. We denote by $\cb$ a
non-degenerate symmetric bilinear form on $\CC^3$ such that
$O(3,\CC)$ is the group of linear transformations of $\CC^3$
preserving this form.
\begin{proposition}\label{P:orbits-flags}
Two complete flags $(E,F)$ and $(E',F')$ in $\CC^3$ belong to the
same $O(3,\CC)$-orbit if and only if the following two conditions
are satisfied:

(i) ranks of the restrictions of $\cb$ to $E$ and to $E'$ are equal;

(ii) ranks of the restrictions of $\cb$ to $F$ and to $F'$ are
equal.
\end{proposition}
\begin{remark}
It is easy to see that the rank of the restriction of $\cb$ to any
line in $\CC^3$ is either 0 or 1, and the rank of the restriction of
$\cb$ to any plane in $\CC^3$ is either 1 or 2.
\end{remark}

Proposition \ref{P:orbits-flags} implies that the $O(3,\CC)$-orbits
on $X$ are described by the following diagram:

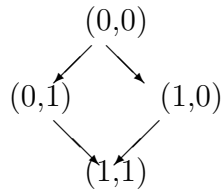
\setlength{\unitlength}{1 cm}
\begin{center}
\begin{picture}(3,3)
      \put(1,0){(1,1)}
      \put(0,1){(0,1)}\put(0.6,0.8){\vector(1,-1){0.6}}
      \put(2,1){(1,0)}\put(2,0.8){\vector(-1,-1){0.6}}
      \put(1,2){(0,0)}\put(1.1,1.8){\vector(-1,-1){0.5}}\put(1.3,1.8){\vector(1,-1){0.5}}
\end{picture}
\end{center}
Let us explain the notation used in this diagram. Here the pair of
numbers $(i,j)$ corresponds to the orbit of flags $(E\subset F)$
such that the restriction of the form $\cb$ to $E$ has a kernel of
dimension $i$, and the restriction of the form $\cb$ to $F$ has a
kernel of dimension $j$. The arrows point from an orbit to another
orbit contained in its closure.

One can easily check the following lemma.
\begin{lemma}\label{4-5}

(1) The orbit $(0,0)$ is open.

(2) The codimensions of the orbits $(0,1)$ and $(1,0)$ are equal to
1.

(3) The codimension of the orbit $(1,1)$ is equal to 2.
\end{lemma}

\hfill

Let us start proving Proposition \ref{P:prop-legth2} by describing
the module $\cm$ more explicitly. Let $U$ be the open
$O(3,\CC)$-orbit in $X$, namely $U$ corresponds to $(0,0)$ in the
above diagram. Let $j:U \hookrightarrow X$ be the identity
imbedding. Fix a flag $(E,F)\in U$. The category of
$O(3,\CC)$-equivariant coherent $\cd$-modules on $U$ is equivalent
to the category of finite dimensional representations of the group
of connected components of the stabilizer of $(E,F)$ in $O(3,\CC)$.
This stabilizer is isomorphic to the group $O(1)\times O(1)\times
O(1)$. The representation of this group corresponding to the
restriction $\cm _0$ of $\cm$ to $U$ is equal to $(A,B,C)\mapsto A$
(recall that $O(1)=\{+1,-1\}$). Then one can easily see that
\begin{eqnarray}\label{M-push}
\cm =j_*\cm_0, \end{eqnarray} where $j_*$ denotes the push-forward
in the category of $\cd$-modules.


\def\co{{\cal O}}


Let us now come back to our situation and let us describe the
perverse sheaf corresponding to our ${\cd}$-module $\cm$. Let us
denote by $M$ (resp. $M_0$) the $O(3,\CC)$-equivariant perverse
sheaf on $X$ (resp. $U$) corresponding to $\cm$ (resp. $\cm_0$) via
the Riemann-Hilbert correspondence (we refer to \cite{borel-book}
for this notion). Clearly $M=j_*M_0$ where $j_*$ is the push-forward
in the derived category of sheaves; since $j\colon U\to X$ is an
affine open imbedding, perverse sheaves go to perverse sheaves (see
\cite{BBD}). It is clear that $M_0$ is a local system on $U$ of rank
one, shifted by 3 to the left. The following lemma describes the
monodromies of $M_0$ around the codimension one orbits $(0,1)$ and
$(1,0)$.
\begin{lemma}\label{L:monodromies}
The monodromy of $M_0$ around the orbit $(1,0)$ is equal to $-1$,
and its monodromy around the orbit $(0,1)$ is equal to $1$.
\end{lemma}

The proof is by a direct computation; it is left to the reader. Let
us now discuss the local system $M_0$ near the minimal orbit
$(1,1)$. Let $\ct$ denote a normal slice at some point on this
orbit. It will be shown below that $\ct$ is isomorphic to $\CC^2$
with the following stratification:
\begin{eqnarray}\label{E:stratification}
\ct=S\cup P\cup Q \cup \{0\},
\end{eqnarray}
where $P$ and $Q$ are parabolas passing through the origin $0$ minus
this origin and such that their closures $\bar P$ and $\bar Q$ are
tangent to each other at 0. $S$ is the open strata which is the
complement of the two parabolas. $P$ is the intersection of $\ct$
with the orbit $(0,1)$, and $Q$ is the intersection of $\ct$ with
the orbit $(1,0)$. Let us describe the stratification
(\ref{E:stratification}) explicitly.

\hfill

Let us choose a basis $e_1,e_2,e_3$ in $\CC^3$ so that the matrix
$\mathcal{B}$ of our quadratic form in this basis is
$$\mathcal{B}=\left[ \begin{array}{ccc}
                0&0&1\\
                0&1&0\\
                1&0&0
                \end{array} \right].$$
Consider the flag $x_0:=(E_0\subset F_0)$ defined
by$E_0=span\{e_1\}, \, F_0=span\{e_1,e_2\}$. Clearly, this flag
belongs to the minimal orbit. By a straightforward computation one
easily checks that the Lie algebra of the stabilizer of $x_0$ in
$O(3,\CC)$ consists of the matrices of the form
$$Lie(St(x_0))=\left\{ \left[ \begin{array}{ccc}
                                  a&b&0\\
                                  0&0&-b\\
                                  0&0&-a
                                  \end{array} \right] \right\}.$$

The transversal in $Lie(O(3,\CC))$ to this subalgebra can be chosen
to be equal to $$\left\{\left[\begin{array}{ccc}
                                  0&0&0\\
                                  -c&0&0\\
                                  0&c&0
                                  \end{array} \right]\right\}.$$
It also can be identified with the tangent space at $x_0$ of the
$O(3,\CC)$-orbit of $x_0$.

Consider now the set of matrices of the form $U=\left[
                                                       \begin{array}{ccc}
                                                       0&0&0\\
                                                       a&0&0\\
                                                       b&a&0
                                                       \end{array}\right]$.
Consider the vectors
$$\xi_1(a,b)=e_1+Ue_1=\left[ \begin{array}{c}
                                                               1\\
                                                               a\\
                                                               b
                                                               \end{array}\right],
                      \xi_2(a,b)=e_2+Ue_2=\left[ \begin{array}{c}
                                                               0\\
                                                               1\\
                                                               a
                                                               \end{array}\right].$$

Consider the flag $X(a,b):=(E(a,b)\subset F(a,b))$, where we denote
$E(a,b)=span\{\xi_1(a,b)\}$ and $F(a,b)=span\{\xi_1(a,b),
\xi_2(a,b)\}$. It is easy to see that the flags $X(a,b)$ with
$a,b\in\CC$, form a transversal to the minimal orbit; this
transversal will be denoted by $\ct$. Let us now compute the
intersection of this transversal with the other orbits. We have:
$$\mathcal{B}(\xi_1(a,b),\xi_1(a,b))=
\left[ \begin{array}{ccc}
             1&a&b \end{array}\right] \left[ \begin{array}{ccc}
                                            0&0&1\\
                                            0&1&0\\
                                            1&0&0
                                            \end{array}\right]
                        \left[ \begin{array}{c}
                                    1\\
                                    a\\
                                    b \end{array}
                                    \right]=2b+a^2,$$

$$\mathcal{B}(\xi_1(a,b),\xi_2(a,b))=
\left[ \begin{array}{ccc}
             1&a&b \end{array}\right] \left[ \begin{array}{ccc}
                                            0&0&1\\
                                            0&1&0\\
                                            1&0&0
                                            \end{array}\right]
                        \left[ \begin{array}{c}
                                    0\\
                                    1\\
                                    a \end{array}
                                    \right]=2a,$$

$$\mathcal{B}(\xi_2(a,b),\xi_2(a,b))=
\left[ \begin{array}{ccc}
             0&1&a \end{array}\right] \left[ \begin{array}{ccc}
                                            0&0&1\\
                                            0&1&0\\
                                            1&0&0
                                            \end{array}\right]
                        \left[ \begin{array}{c}
                                    0\\
                                    1\\
                                    a \end{array}
                                    \right]=1.$$
 Thus we obtain that the matrix of the
 restriction of the quadratic form $B$ to $F(a,b)$ is equal to
 $\left[ \begin{array}{cc}
             2b+a^2&2a\\
             2a&1 \end{array}\right]$.

Thus we see that the intersection of the transversal with the
closure of the orbit $(0,1)$ consists of pairs $\bar
P:=\{(a,b)|2b=3a^2\}$, and the intersection of the transversal with
the closure of the orbit $(1,0)$ consists of the pairs $\bar
Q:=\{(a,b)|2b=-a^2\}$. In particular we see that $\bar P$ and $\bar
Q$ are two different parabolas on $\CC^2$ with a single common point
$(0,0)$ and which are tangent to each other at this point. To be
consistent with the above notation, denote $P:=\bar P\backslash
\{0\}$, $Q:=\bar Q\backslash \{0\}$.

Let us agree on the following notation: for any
$O(3,\CC)$-equivariant perverse sheaf $T$ on $X$ or on a locally
closed subset of $X$, we will denote by $\tilde T$ the restriction
of $T$ to the intersection of the transversal $\ct$ with this
subset. Then Lemma \ref{L:monodromies} implies that in this notation
$\tilde M_0$ is a rank one local system on $S$ (up to a
cohomological shift) with monodromy $+1$ around $P$, and monodromy
$-1$ around $Q$. Let $\tilde j\colon S\to \ct$ denote the identity
imbedding. It is easy to see that the restriction
$\widetilde{j_*M_0}$ is equal to $\tilde j_*\tilde M_0$. We will
need the following lemma.
\begin{lemma}\label{L:length}
The perverse sheaf $\tilde j_*\tilde M_0$ has length 3. Moreover the
smooth parts of all the irreducible subquotients are supported on
different strata, and they are rank one local systems (up to a
cohomological shift).
\end{lemma}
Before we prove Lemma \ref{L:length}, let us finish the proof of
Theorem \ref{4-3}. It is easy to see that the functor of restriction
to the transversal $\ct$ is an exact and faithful functor from the
category of $O(3,\CC)$-equivariant perverse sheaves on $X$ to the
category of perverse sheaves on $\ct$ (this is because $\ct$
intersects every orbit). This and Lemma \ref{L:length} imply that
$j_*M_0$ has length at most 3 as $O(3,\CC)$-equivariant perverse
sheaf. Hence it remains to show that the restriction to $\ct$ of any
irreducible subquotient of $j_*M_0$ (as $O(3,\CC)$-equivariant
perverse sheaf) is an irreducible perverse sheaf; that clearly will
imply that the length of $j_*M_0$ is equal to 3.

Let $N$ be an irreducible subquotient of $j_*M_0$ as an
$O(3,\CC)$-equivariant perverse sheaf. Consider its restriction
$\tilde N$ to $\ct$. We want to prove that $\tilde N$ is
irreducible. Assume in the contrary that it is not irreducible. By
Lemma \ref{L:length} the smooth parts of different irreducible
subquotients of $\tilde N$ are supported on different strata, and
all of these smooth parts are rank one local systems.

Any irreducible $O(3,\CC)$-equivariant perverse sheaf, can be
constructed as follows. There exists a unique $O(3,\CC)$-orbit
$\co\subset X$ and an irreducible $O(3,\CC)$-equivariant local
system $L$ on $\co$ such that $N$ is isomorphic to the
Goresky-MacPherson extension of $L$. It is easy to see that $\tilde
N$ is the Goresky-MacPherson extension of $\tilde L$.

Let us show that $rk (L)=1$. Assume in the contrary that $rk(L)\geq
2$. Then $\tilde N$ is not irreducible by Lemma \ref{L:length},
moreover there are at least $rk(L)$ irreducible subquotients such
that their smooth part is a rank one local system supported on
$\co\cap \ct$. This contradicts Lemma \ref{L:length}, according to
which the smooth parts of all irreducible subquotients of $\tilde N$
must be supported on different orbits.

Thus $\tilde L$ is a rank one local system. Its Goresky-MacPherson
extension must be irreducible. Thus $\tilde N$ is irreducible. Thus
we have shown that the restriction to $\ct$ maps irreducible
$O(3,\CC)$-equivariant subquotients of $j_*M_0$ to irreducible
perverse sheaves. Hence length of $M=j_*M_0$ is equal to 3.

Now it remains to prove Lemma \ref{L:length}.

\hfill

{\bf Proof of Lemma \ref{L:length}.} It will be more convenient to
make a change of variables on the transversal $\ct$ as follows:
\begin{eqnarray*}
z=a,\\
w=\frac{1}{2}(\frac{a^2}{2}+b).
\end{eqnarray*}
In coordinates $(z,w)$
$$\bar Q=\{w=0\}, \bar P=\{w=z^2\}.$$
Recall that on $S=\ct\backslash(\bar P\cup \bar Q)$ we have a rank
one local system $\tilde M_0$ whose monodromy around $\bar P$ is
equal to $+1$, and around $\bar Q$ to $-1$. Let $j\colon S\to \ct$
be the identity imbedding. We want to show that $j_*\tilde M_0$ is a
length 3 perverse sheaf, and the smooth parts of different
irreducible subquotients are (appropriately shifted) rank one local
systems with different supports.

Let $f\colon S\to \ct\backslash \bar Q$ be the identity imbedding,
which is an affine morphism. Since the monodromy of $\tilde M_0$
around $P$ is equal to $+1$, there is a short exact sequence of
perverse sheaves:
\begin{eqnarray}\label{E:short-exact}
K[2]\to f_*\tilde M_0\to N[1],
\end{eqnarray}
where $K$ is a rank one local system on $\ct\backslash \bar Q$, and
$N$ is a rank one local system on $P(=\bar P\backslash \{0\})$.

Let $g\colon\ct\backslash \bar Q\to \ct$ be the identity imbedding.
It is an affine imbedding, hence $g_*$ is an exact functor between
categories of perverse sheaves. Also $j=g\circ f$. Hence applying
$g_*$ to (\ref{E:short-exact}) we get a short exact sequence of
perverse sheaves on $\ct$
\begin{eqnarray}\label{E:short-exact-2}
g_*K[2]\to j_*\tilde M_0[2]\to j_*N[1].
\end{eqnarray}
It is easy to see that the monodromy of $K$ around $\bar Q$ is equal
to $-1$. Hence $g_*K[2]$ is an irreducible perverse sheaf; clearly its
smooth part is the rank one local system $K[2]$ supported on
$\ct\backslash \bar Q$.

It remains to show that $j_*N[1]$ has length 2, and the smooth parts of
the two irreducible subquotients are (shifted) rank one local
systems on $P$ and $\{0\}$ respectively. Let us denote by $h\colon
P\to \bar P$ the identity map. It is equivalent to prove that
$h_*N[1]$ has length 2, and the irreducible subquotients have
appropriate supports. It suffices to show that the monodromy of $N$
around $0\in\bar P$ is equal to $+1$. To prove this, let us fix the
loop around 0 in $P$ given by
\begin{eqnarray}\label{E:loop1}
\{(e^{i\theta},e^{2i\theta})|\, 0\leq \theta\leq 2\pi\}.
\end{eqnarray}
Around each point of this loop let us fix a small loop
\begin{eqnarray}\label{E:loop2}
\{(e^{i\theta},e^{2i\theta}+\eps e^{i\omega})|\, 0\leq \omega\leq
2\pi\},
\end{eqnarray}
where $0<\eps<1$ is fixed, $\theta\in [0,2\pi]$ is arbitrary. It is
easy to see that each point on the loops (\ref{E:loop2}) belongs to
the open stratum $S\subset \ct$, and all these points together form
a 2-dimensional real torus $T^2$. For a fixed $\theta$, the loop
(\ref{E:loop2}) is a loop in $S$ around $P$.

Let us write the torus $T^2$ as a product of circles $T^2=S_1\times
S_2$, where $S_1$ corresponds to the parameter $\theta$ in
(\ref{E:loop2}), and $S_2$ to $\omega$. Let $\pi_1\colon T^2\to S_1$
denote the natural projection.

Let us denote by $\bar M_0$ the restriction of the (shifted by 2)
local system $M_0$ to the torus $T^2$. Clearly the restriction of
$N[1]$ to the circle (\ref{E:loop1}) is equal to $R^1\pi_{1*}(\bar
M_0)$ (where $R^1\pi_{1*}$ is the first right derived functor of the
push-forward under $\pi_1$) which is a rank one local system on the
circle $S_1$. We want to show that its monodromy is equal to $+1$.

Any local system on torus $T^2$ is characterized (up to isomorphism)
by its monodromies along the meridians. The monodromy of $\bar M_0$
along the meridian $S_2$ is equal to $+1$. Then evidently the
monodromy of $R^1\pi_{1*}\bar M_0$ is equal to the monodromy of
$\bar M_0$ along the meridian $S_1$, which we can choose to be
\begin{eqnarray}\label{E:loop3}
\{(e^{i\theta},e^{2i\theta}+\eps)|\, \theta\in [0,2\pi]\}.
\end{eqnarray}
Thus it remains to show that the monodromy of $\bar M_0$ along the
circle (\ref{E:loop3}) is equal to $+1$. To do that, let us go back
from the torus to $\CC^2=\ct$. In $\CC^2$ we can close the loop
(\ref{E:loop3}) by a disk $D$ defined by
$$D:=\{(z,z^2+\eps)|\, |z|\leq 1\}.$$
Clearly $\pt D$ equals the circle (\ref{E:loop3}); $D$ does not
intersect $\bar P$; and does intersect $\bar Q$ in exactly two
points $(\pm i \sqrt{\eps},0)$, in both of them the intersection is
transversal. The latter implies that the monodromy of $\bar M_0$
around $\pt D$ is equal to the product of monodromies of the two
loops, each surrounding once exactly one of the the two points of
intersection. But the monodromy of each such loop is equal to $-1$,
hence their product is $+1$. Hence the monodromy of $\bar M_0$
around (\ref{E:loop3}) is $+1$. Hence the monodromy of $N$ around
$0\in \bar P$ is $+1$. Lemma \ref{L:length} is proved.




\end{document}